\documentclass[english,10pt,a4paper]{article}
\usepackage[utf8]{inputenc}
\usepackage[T1]{fontenc}
\usepackage{babel}
\usepackage{hyperref}

\usepackage{textcomp}
\usepackage{graphicx}%
\usepackage{multirow}%
\usepackage{amsmath,amssymb,amsfonts}%
\usepackage{amsthm}%
\usepackage{mathrsfs}%
\usepackage[title]{appendix}%
\usepackage{xcolor}%
\usepackage{textcomp}%
\usepackage{manyfoot}%
\usepackage{booktabs}%
\usepackage{algorithm}%
\usepackage{algorithmicx}%
\usepackage{algpseudocode}%
\usepackage{listings}%
\theoremstyle{definition}
\newtheorem{Theorem}{Theorem}[section]
\newtheorem{Definition}{Definition}[section]
\newtheorem{Lemma}{Lemma}[section]
\newtheorem*{Proof}{Proof}

\usepackage{tikz-cd}
\usepackage{authblk}
\usepackage[numbers,sort&compress]{natbib}
\begin{document}
	\title{Compatibility axioms for left-regular bands to construct Hopf algebras}
	\author{Lingxiao Hao\thanks{Corresponding author: 21110180006@m.fudan.edu.cn}}
	
	\author{Shenglin Zhu}
	\affil{School of Mathematical Sciences, Fudan University, 220 Handan Road, Shanghai, 200433, China} 

	\date{\today}
	\maketitle
	\begin{abstract}
		Aguiar and Mahajan provided several coalgebra axioms and algebra axioms for a family of left-regular bands to construct a commutative diagram of algebras and coalgebras. In this paper, we will add compatibility axioms to make it a diagram of Hopf algebras.
	\end{abstract}
	\section{Introduction}
	Left regular bands, or LRBs for short, are semigroups with relations $ x^{2}=x,$ $xyx=xy$. The
	origin of LRBs can be traced to Klein-Barmen \cite{ref11} and Sch\"{u}tzenberger \cite{ref6}. A significant example of a LRB is the Coxeter complex of a Coxeter group.\par
	Aguiar and Mahajan \cite{ref1} offered several coalgebra axioms and algebra axioms for a family of left-regular bands $\left\lbrace \Sigma^{n}  \right\rbrace_{n\geq 0}$ to construct a diagram of algebras and coalgebras: 	\begin{equation}
			\begin{tikzcd}\label{vs}
			\mathop{\oplus}\limits_{n\geqslant 0}\mathbb{K}\Sigma^{n} \arrow[r, "base^{\ast}"] \arrow[dd, "supp"', two heads] \arrow[rrdd, "\Upsilon^{\ast}"'] & \mathop{\oplus}\limits_{n\geqslant 0}\mathbb{K}Q^{n} \arrow[r, "\Theta"] \arrow[rrdd, "\Psi"] \arrow[dd, "lune", two heads] & \mathop{\oplus}\limits_{n\geqslant 0}\mathbb{K} \left(\mathcal{C}^{n}\times\mathcal{C}^{n}\right) \arrow[rd, "s"] &                            \\
			&                                                                                      &                            & \mathop{\oplus}\limits_{n\geqslant 0}\mathbb{K}\left( \mathcal{C}^{n}\times\mathcal{C}^{n}\right) ^{\ast} \arrow[d, "Road"]        \\
			\mathop{\oplus}\limits_{n\geqslant 0}\mathbb{K} L^{n} \arrow[r, "base^{\ast}"] \arrow[rrd, "\Phi"']                              & \mathop{\oplus}\limits_{n\geqslant 0}\mathbb{K} Z^{n} \arrow[rrd, "\Upsilon" description]                       & \mathop{\oplus}\limits_{n\geqslant 0}\mathbb{K} (Z^{n})^{\ast} \arrow[r, "lune^{\ast}" description, hook] \arrow[d, "base" description]                              & \mathop{\oplus}\limits_{n\geqslant 0}\mathbb{K}\left( Q^{n}\right) ^{\ast} \arrow[d, "base"] \\&
			& \mathop{\oplus}\limits_{n\geqslant 0}\mathbb{K} (L^{n})^{\ast} \arrow[r, "supp^{\ast}"', hook]                               & \mathop{\oplus}\limits_{n\geqslant 0}\mathbb{K}\left( \Sigma^{n}\right) ^{\ast}               
		\end{tikzcd}
	\end{equation}
	where the notations will be explained in the third section.
	In particular, when $\Sigma^{n}$ is the set of compositions of $[n]$, which can be seen as the Coxeter complex of the Coxeter group of type $A_{n-1}$ (that is, the symmetry group $S_{n}$), the above diagram projects onto the following diagram \cite[Theorem 6.1.4]{ref1}:
	\begin{equation*}
			\begin{tikzcd}
			N\Lambda \arrow[r, hook] \arrow[d, two heads] & S\Lambda \arrow[d, two heads] \\
			\Lambda \arrow[r, hook]                       & Q\Lambda                     
		\end{tikzcd}
	\end{equation*}
where $S\Lambda$ is the Malvenuto-Reutenauer Hopf algebra \cite{ref10}, $\Lambda$, $Q\Lambda$ and $N\Lambda$ are the Hopf algebras of symmetric functions \cite{ref7},  quasi-symmetric functions \cite{ref36} and non-commutative symmetric functions \cite{ref36} respectively.
\par Aguiar and Mahajan raised a question \cite[Section 6.1]{ref1} : ``For a family $\left\lbrace \Sigma^{n}  \right\rbrace_{n\geq 0}$ of LRBs, can one give compatibility axioms
such that if a family satisfies the coalgebra, algebra and compatibility axioms then Diagram \ref{taget} is a diagram of Hopf algebras?'' Our work is to solve this problem.
\par In the second section, we recall some basic definitions related to Hopf algebras, LRBs and posets. In the third section, we will introduce Diagram \ref{taget} and  several coalgebra and algebra axioms provided by Aguilar and Mahajan. 
\par	The fourth section presents our main results. We will provide compatibility axioms and prove they can make \ref{taget} a diagram of connected graded Hopf algebras. We will give an example to prove the compatibility axioms are  compatible with the algebra and coalgebra axioms.
	\section{Preliminaries}

\par We begin by recalling some basic facts about Hopf algebras. Readers are referred to the books of Montgomery \cite{Montgomery} and Sweedler \cite{ref3} for more details.
   
	Let $k$ be a field and $\circ$ be the composition of maps.
		A \emph{$k$-coalgebra} is a $k$-vector space $C$ together with two $k$-linear maps, \emph{comultiplication} $\bigtriangleup : C\rightarrow C \otimes C $ and \emph{counit} 
		$\varepsilon: C\rightarrow k$, such that the following diagrams are commutative:
		\begin{equation*}
		\begin{tikzcd}
			C \arrow[d, "\bigtriangleup"] \arrow[r, "\bigtriangleup"] & C\otimes C \arrow[d, "id \otimes\bigtriangleup"] & k\otimes C & C \arrow[l, "1\otimes"'] \arrow[r, "\otimes 1"] \arrow[d, "\bigtriangleup"']        & C\otimes k \\
			C\otimes C \arrow[r, "\bigtriangleup\otimes  id"]         & C\otimes C\otimes C                              &            & C\otimes C \arrow[lu, "\varepsilon\otimes id"] \arrow[ru, "id \otimes\varepsilon"'] &           
		\end{tikzcd}		\end{equation*}
	
	For coalgebras $C$ and $D$, with comultiplications $\bigtriangleup_{C}$
	and $\bigtriangleup_{D}$, and counits $\varepsilon_{C}$ and $\varepsilon_{D}$ respectively, a map $f: C\rightarrow D$  is a \emph{coalgebra morphism} if $\bigtriangleup_{D}\circ f = (f \otimes f)\circ\bigtriangleup_{C}$
	and if $\varepsilon_{C} =\varepsilon_{D} \circ f$.
A subspace $I\subseteq C$ is a \emph{coideal} if $\bigtriangleup I\subseteq I\otimes C + C \otimes I$ and if $\varepsilon(I) = 0$.	
	It is easy to check that if $I$ is a coideal, then the $k$-space $C / I$ is a coalgebra
	with comultiplication induced from $\bigtriangleup$, and conversely.\par
	For a coalgebra $C$, $C\otimes C$ is a coalgebra with $$\bigtriangleup _{C\otimes C}:=(id\otimes\tau\otimes id)\circ(\bigtriangleup\otimes\bigtriangleup)$$where $\tau:C \otimes C\rightarrow C\otimes C$ is the \emph{twist} map given by $\tau(v\otimes w) = w \otimes v$, and $\varepsilon_{C\otimes C}$ given by $$\varepsilon_{C\otimes C}(a\otimes b)=\varepsilon(a)\varepsilon(b) .$$	\par The field k is a coalgebra with the coproduct $\bigtriangleup_{k}$ given by $1\mapsto1\otimes 1$ and $\varepsilon_{k}$ given by $1\mapsto 1.$\par
	A $k$-space $B$ is a \emph{bialgebra} if $(B, m, u)$ is an algebra ($m$ is the product and $u:k\rightarrow B$ is the unit), $ (B,\bigtriangleup,\varepsilon) $ is a coalgebra, and $m$ and $u$ are coalgebra morphisms.
			A map  $f : B\rightarrow B'$  between bialgebras is called a \emph{bialgebra morphism} if it is both an algebra morphism and a coalgebra morphism, and a subspace 
		$I \subseteq B$ is a \emph{biideal} if it is both an ideal and a coideal. 
	The quotient $B / I$ is a 
	bialgebra precisely when $I$ is a biideal of $B$.
	\par We will need the following definition for later proof.
	\begin{Definition}
		Let $B$ be a bialgebra. $B$ is \emph{connected graded} if $B=\oplus_{n\geq 0}B_{n}$, $B_{0}=k$, $B_{i}\cdot B_{j}\subseteq B_{i+j}$, and $\bigtriangleup B_{n}\subseteq \bigoplus_{i+j=n}B_{i}\otimes B_{j}.$ 
	\end{Definition}
	
For a bialgebra	$(H, m, u,\bigtriangleup,\varepsilon) $, $H$ is a \emph{Hopf algebra} if there exists an element $S \in Hom_{k}(H, H)$ s.t. $$m\circ(S\otimes id)\circ\bigtriangleup=u\circ\varepsilon.$$ $S$ is called an \emph{antipode} for $H$.
For Hopf algebras $H, K$ be  with antipodes $S_{H}$, $S_{K}$ respectively, a map 
 $f : H\rightarrow K$ is a \emph{Hopf morphism} if it is a bialgebra morphism 
and $f\circ S_{H}=S_{K}\circ f$.

We also need to recall some notions related to LRBs. A \emph{left-regular band} (LRB for short)\cite{ref6, ref11}, is a semigroup $\Sigma$ which satisfies the identities$$ x^{2}=x, xyx=xy $$ for all $ x, y \in \Sigma $. \par
$\Sigma$ is a poset with the partial order \cite[Section 2.2.1]{ref1} defined by $$x\leq y\Leftrightarrow xy=y.$$\par Elements of $\Sigma$  are called \emph{faces} and we say that $x$ is a face of $y$ if $x\leq y$.\par
Maximal elements in $\Sigma$ are called \emph{chambers}. The set of chambers is denoted by $\mathcal{C}$.\par
Another relation $\preceq$ \cite[Section 2.2.2]{ref1} on $\Sigma$ is defined by $x\preceq y \Leftrightarrow yx = y. $ It is
transitive and reflexive, but not necessarily anti-symmetric. 
This gives an equivalent relation $\sim$ by $$x\sim y \Leftrightarrow x\preceq y, y\preceq x. $$
Let $L=\Sigma/\sim$. Denote the quotient map by $ \emph{supp}:\Sigma \twoheadrightarrow L$. Then one has
$$	yx = y \Leftrightarrow supp(x) \leq supp( y).$$The map ``$supp$'' is order-preserving.\par

\par Elements of $L$ are called \emph{flats}. 
There is an analogue of the last subsection with faces and flats replaced
by pointed faces and lunes respectively.	Let  $Q= \left\lbrace (x, c) \in\Sigma\times\mathcal{C}| x \leq c\right\rbrace. $ \par There is a partial order on $Q$ defined by
$(x, c) \leq (y, d) \Leftrightarrow c = d$  and $ x \leq y.$\par
Elements of $Q$ are called \emph{pointed faces}.
Another relation $\preceq$ on $Q$ is given by $(x, c)\preceq (y, d) \Leftrightarrow yx =
y$ and $yc = d$. It is transitive and reflexive, but not necessarily anti-symmetric. This relation gives an equivalent relation $\sim$ by $$(x, c)\sim(y, d) \Leftrightarrow (x, c)\preceq (y, d), (y, d)\preceq (x, c). $$
Let $Z=Q/\sim$. The quotient map is denoted by \emph{lune}: $Q\twoheadrightarrow Z$. Then
$$yx = y, yc = d\Leftrightarrow lune(x, c) \leq lune(y, d)$$
holds by definition.\par Elements of $Z$ are called \emph{lunes} and $lune(x, c)$ is called the lune	of $x$ and $c$. The $lune$ map preserves order.\par The map
$reg:Q \rightarrow\left\lbrace  R | R \subseteq \Sigma\right\rbrace $ is given by
$reg(x, c) = \left\lbrace y | xy \leq c\right\rbrace$. Then
\cite[Lemma 2.3.1]{ref1} there is a commutative diagram\\
\begin{equation*}
	\begin{tikzcd}
		& Q \arrow[ld, "lune", two heads] \arrow[rd, "reg", two heads] &    \\
		Z \arrow[rr, "zone", two heads] &                                                              & Z'.
	\end{tikzcd}
\end{equation*}
This diagram gives the definition of the zone map.
The following lemma will be useful in later proof.
\begin{Lemma}\cite[Lemma 2.7.4]{ref1}
	If $\Sigma$ is a LRB with a minimum element $\emptyset$, then $\emptyset$ is the unit of the product.
\end{Lemma}
In addition, some basic notions related to graded poset are needed as well. Let $P$ be a poset. An element $y$ in $P$ \emph{covers} another element $x$ in $P$
if $x<y$ and there is no $z\in P$ such that $x<z<y$. $P$ is \emph{ranked} if there exists a map $\rho:P\rightarrow \mathbb{Z}$ s.t. $\rho(y)=\rho(x)+1$ if $y$ covers $x$. The \emph{rank} of an element $x$ is $\rho(x)$. A ranked poset with minimum element $\hat{0}$ is \emph{graded} if $\rho(\hat{0})=0.$ A graded poset $P$ is \emph{of rank n} if the rank of every maximal element in $P$ is $n$.

	\section{The problem}
	
\subsection{Coalgebra and algebra axioms}	
In this subsection, we will introduce the coalgebra and algebra axioms given by \cite[Chapter 6]{ref1}. \par
Let $ \left\lbrace \Sigma^{n}  \right\rbrace_{n\geq 0} $ be a family of LRBs, where $\Sigma^{n}$ is a finite graded poset of rank $\max\left\lbrace n-1,0 \right\rbrace $ and has a unique minimum element denoted by $ \emptyset_{n}$. 
 Let $\Sigma^{1}=\left\lbrace \emptyset_{1} \right\rbrace$ and $\Sigma^{0}=\left\lbrace \emptyset_{0} \right\rbrace $. Let $\mathcal{C}^{n}$ be the set of maximal elements in $\Sigma^{n}$. For an element $K\in \Sigma^{n}$, let $deg K := rank K + 1,$ where ``$rank$'' denotes the rank of an element. Let $\Sigma^{n}_{K}:=\left\lbrace F\in \Sigma^{n}|K\leq F \right\rbrace$, $  \mathcal{C}^{n}_{K}:=\left\lbrace F\in \mathcal{C}^{n}|K\leq F \right\rbrace. $	
The following axioms were given in \cite[Chapter 6]{ref1}. $(C1)-(CP)$ \cite[6.3.1]{ref1}  are called the coalgebra axioms while $(A1)-(AP)$ \cite[6.6.1]{ref1} are called the algebra axioms. We modify them slightly for trivial cases. Our modifications won't influence the results in \cite[Chapter 6]{ref1}.\\ \hspace*{\fill} \\
\textbf{Axiom }$(C1).$  For every $F\in \Sigma^{n}$, there exists a composition $(f_{1},f_{2},\cdots,f_{degF})$ of $n$ and a poset isomorphism$$b_{F}:\Sigma^{n}_{F}\rightarrow \Sigma^{f_{1}}\times \Sigma^{f_{2}}\times \cdots\times\Sigma^{f_{degF}}.$$\par For $F=\emptyset_{n}$, the composition is $(n)$ and the poset isomorphism $b_{\emptyset_{n}}$ is the identity map id$:=\Sigma^{n}\rightarrow\Sigma^{n}$.\\ \hspace*{\fill} \\
\textbf{Axiom }$(C2)$.  
Let $K \leq F$. First apply Axiom $(C1)$ to  $K \in \Sigma^{n}$ to get a composition $(k_{1},k_{2},\cdots,k_{degK})$ of $n$ and a poset isomorphism
$$b_{K}:\Sigma^{n}_{K}\rightarrow \Sigma^{k_{1}}\times \Sigma^{k_{2}}\times \cdots\times\Sigma^{k_{degK}},$$
which maps $F$ to $F_{1} \times F_{2}\times \cdots \times F_{deg K}.$ This induces an isomorphism $$b_{K}:\Sigma^{n}_{F}\rightarrow \Sigma^{k_{1}}_{F_{1}}\times \Sigma^{k_{2}}_{F_{2}}\times \cdots\times\Sigma^{k_{degK}}_{F_{degK}}.$$ \par    
Apply Axiom $(C1)$ to $F_{i}\in\Sigma^{k_{i}} $ to get a composition $(f_{i1},f_{i2},\cdots,f_{idegF_{i}})$
of $ k_{i} $ and a poset isomorphism$$b_{F_{i}}:\Sigma^{k_{i}}_{F_{i}}\rightarrow \Sigma^{f_{i1}}\times \Sigma^{f_{i2}}\times \cdots\times\Sigma^{f_{idegF_{i}}},$$ for $1 \leq i \leq deg K$.
Then, we have$$b_{F} = (b_{F_{1}} \times b_{F_{2}}\times \cdots \times b_{F_{deg K}} ) \circ b_{K},$$
where $b_{F}$ is the poset map given by Axiom $(C1)$ for $F$. In particular,
we require that
$$(f_{1}, f_{2},\cdots,f_{deg F} )=(f_{11}, f_{12},\cdots,f_{1deg F_{1}}, f_{21},\cdots,f_{mdeg F_{m}} ),$$ where $m = deg K$.
\\ \hspace*{\fill} \\
\textbf{Axiom }$(C3).$  Every rank 2 face $H$ of $\Sigma^{n}$ has exactly two rank 1 faces, say $K$ and $G$. Further the poset isomorphisms $b_{K} : \Sigma^{n}_{K}\rightarrow \Sigma^{k_{1}} \times\Sigma^{k_{2}} $and  $b_{G} : \Sigma^{n}_{G}\rightarrow \Sigma^{g_{1}} \times\Sigma^{g_{2}} $
map $H$ to $\emptyset_{k_{1}}\times K'$ and $G'\times\emptyset_{g_{2}}$ respectively, for some rank 1 faces $K' \in \Sigma^{k_{2}}$ and
$G' \in \Sigma^{g_{1}}$.
	\\ \hspace*{\fill} \\
	\textbf{Axiom }$(CP).$     	For
	any $F \leq H, N \in \Sigma^{n},$ we have
	$$ b_{F} (HN) = b_{F} (H)b_{F} (N),$$
	where $b_{F}$ is the poset isomorphism associated to $F$ by Axiom $(C1)$.\\ \hspace*{\fill} \\
\textbf{Axiom }$(A1).$ For every composition $(f_{1}, f_{2},\cdots,f_{k})$ of $n$, there exists a face $F$ of
$\Sigma^{n}$ of degree $k$, and a poset isomorphism
$$j_{F} : \Sigma^{f_{1}}\times \Sigma^{f_{2}}\times \cdots\times\Sigma^{f_{degF}} \rightarrow \Sigma^{n}
_{F} .$$\par
Further, distinct compositions give distinct faces; so it is unambiguous to use the
notation $j_{F}$ for the poset map. If the composition is $(n)$ then $F =\emptyset_{n}$ and the poset
isomorphism $j_{\emptyset_{n}}$ is the identity map id : $\Sigma^{n} \rightarrow \Sigma^{n}.$\par
In addition, for a sequence $(g_{1}, g_{2},\cdots,g_{m})$ of non-negative integers, identify $(g_{1}, g_{2},\cdots,g_{m})$ with $(f_{1}, f_{2},\cdots,f_{k})$ as a composition, where $(f_{1}, f_{2},\cdots,f_{k})$ is given by deleting all the zeros in $(g_{1}, g_{2},\cdots,g_{m})$. The composition of the projection $\Sigma^{g_{1}}\times \Sigma^{g_{2}}\times \cdots\times\Sigma^{g_{m}} \rightarrow\Sigma^{f_{1}}\times \Sigma^{f_{2}}\times \cdots\times\Sigma^{f_{degF}}  $ and $j_{F} $ is still denoted $j_{F} $. \\ \hspace*{\fill} \\
\textbf{Axiom }$(A2).$ 
Let $(g_{1}, g_{2},\cdots,g_{m})$ be a composition of $n$ and for $1\leq i\leq m$, let $(f_{i1}, f_{i2},\cdots,f_{ik_{i}} )$
be a composition of $g_{i}$. First apply Axiom $(A1)$ to the composition of $g_{i}$ to get
a face $F_{i}$ of $\Sigma^{n}$ of degree $k_{i}$, and a poset isomorphism
$$j_{F_{i}} : \Sigma^{f_{i1}}\times \Sigma^{f_{i2}}\times \cdots\times\Sigma^{f_{ik_{i}}}\rightarrow \Sigma ^{g_{i}}
_{ F_{i}}$$ for $1 \leq i \leq m$.
Apply Axiom $(A1)$ to the composition $(g_{1}, g_{2},\cdots,g_{m})$ of $n$ to get a face $G$
of degree $m$ and a poset isomorphism
$$j_{G} : \Sigma^{g_{1}}\times \Sigma^{g_{2}}\times \cdots\times\Sigma^{g_{m}} \rightarrow \Sigma^{n}
_{G},$$
which maps $F_{1} \times F_{2 }\times \cdots \times F_{m}$ to $F'.$ This induces an isomorphism$$j_{G} : \Sigma^{g_{1}}_{F_{1}}\times \Sigma^{g_{2}}_{F_{2}}\times \cdots\times\Sigma^{g_{m}}_{F_{m}} \rightarrow \Sigma^{n}
_{F'},$$
with $G \leq F'$ .
Consider the composition of $n$
$$(f_{11}, f_{12},\cdots,f_{1k_{1}}, f_{21},\cdots,f_{2k_{2}},\cdots,f_{m1},\cdots,f_{mk_{m}}),$$
which is the concatenation of the $m$ compositions that we started with. Apply the
axiom $(A1)$ to it to get a face $F \in \Sigma^{n}$. Then
$j_{F} = j_{G} \circ (j_{F_{1}}\times j_{F_{2}} \times \cdots \times j_{F_{m}}).$ 
In particular, we require that $F = F'$ .
\\ \hspace*{\fill} \\
\textbf{Axiom }$(AP).$ 
$$j_{F}
^{-1}(HN) = j_{F}
^{-1}(H)j_{F}
^{-1}(N),$$
for $F,H, N\in \Sigma^{n}$, $F\leq H,N$.
\\ \hspace*{\fill} \\
\subsection{The diagram}

Let $\mathbb{K} $ be a field of characteristic 0. In this subsection, we will give the vector spaces in \ref{vs} structures of algebras and coalgebras. All of these definitions are from \cite[Chapter 6]{ref1}. \par	Let $ \left\lbrace \Sigma^{n}  \right\rbrace_{n\geq 0} $ be a family of LRBs, which satisfies all coalgebra axioms $(C1)-(CP)$ and all algebra axioms $(A1)-(AP)$.\par Let $Q^{n}=\left\lbrace (F,D)\in \Sigma^{n}\times\mathcal{C}^{n}|F\leq D \right\rbrace, $ 
for $n\geq 0$. Let $L^{n}=supp(\Sigma^{n})$, $Z^{n}=lune(Q^{n})$, for $n\geq 0.$
\begin{itemize}
	\item[$\bullet$]\cite[Sections 6.5.3 and 6.8.2]{ref1} Let $\mathcal{P}:=\mathop{\oplus}\limits_{n\geqslant 0}\mathbb{K}\left( \Sigma^{n}\right) ^{\ast}$, with the basis element corresponding to
	$F\in\Sigma^{n}$ denoted by $M_{F}$ when $n>0$ and 1 when $F=\emptyset_{0}$.\par
	 The coproduct on $\mathcal{P} $ is given by$$\bigtriangleup\left(M_{F}\right) =1\otimes M_{F}+M_{F}\otimes 1+\mathop{\sum}\limits_{K:rankK=1, K\leq F}M_{F_{1}}\otimes M_{F_{2}}, $$where $b_{K}(F)=F_{1}\times F_{2}. $\par
	The product on $\mathcal{P} $ is given by
	$$M_{F_{1}}\ast M_{F_{2}}= \mathop{\sum}\limits_{F: GF=j_{G}(F_{1}\times F_{2})}M_{F}.$$The vertex $G$ is  given by
	Axiom $(A1)$ for the composition $(g_{1}, g_{2})$, where $F_{i} \in \Sigma^{g_{i}}$.
	\item[$\bullet$] \cite[Sections 6.5.5 and 6.8.4]{ref1}Let $\mathcal{Q}:=\mathop{\oplus}\limits_{n\geqslant 0}\mathbb{K}\left( Q^{n}\right) ^{\ast}$, with the basis element corresponding to
	$(F,D)\in Q^{n}$ denoted by $M_{(F,D)}$ when $n>0$ and 1 when $(F,D)=(\emptyset_{0},\emptyset_{0})$.\par
The coproduct on $\mathcal{Q} $ is given by$$\bigtriangleup\left(M_{(F,D)}\right) =1\otimes M_{(F,D)}+M_{(F,D)}\otimes 1+\mathop{\sum}\limits_{K:rankK=1, K\leq F}M_{(F_{1},D_{1})}\otimes M_{(F_{2},D_{2})}, $$where $b_{K}(F)=F_{1}\times F_{2}$, $ b_{K}(D)=D_{1}\times D_{2}. $\par
	The product on $\mathcal{Q} $ is given by
	$$M_{(F_{1},D_{1})}\ast M_{(F_{2},D_{2})}= \mathop{\sum}\limits_{F: GF=j_{G}(F_{1}\times F_{2})}M_{(F,Fj_{G}(F_{1}\times F_{2}))}.$$The vertex $G$ is  given by
	Axiom $(A1)$ for the composition $(g_{1}, g_{2})$, where $F_{i} \in \Sigma^{g_{i}}$, $i=1,2$.
	\item[$\bullet$] \cite[Sections 6.5.4 and 6.8.3]{ref1}
	Let $\mathcal{M}:=\mathop{\oplus}\limits_{n\geqslant 0}\mathbb{K}\Sigma^{n} $, with the basis element corresponding to
	$P\in\Sigma^{n}$ denoted by $H_{P}$ when $n>0$ and 1 when $P=\emptyset_{0}$.\par
	The coproduct on $\mathcal{M} $ is given by$$\bigtriangleup\left(H_{P}\right) =1\otimes H_{P}+H_{P}\otimes 1+\mathop{\sum}\limits_{K:rankK=1}H_{P_{1}}\otimes H_{P_{2}}, $$where $b_{K}(KP)=P_{1}\times P_{2}. $\par
	The product on $\mathcal{M} $ is given by
	$$H_{P_{1}}\ast H_{P_{2}}=H_{j_{G}(P_{1}\times P_{2})}.$$The vertex $G$ is  given by
	Axiom $(A1)$ for the composition $(g_{1}, g_{2})$, where $P_{i} \in \Sigma^{g_{i}}$, $i=1,2$.
	\item[$\bullet$]\cite[Sections 6.5.6 and 6.8.5]{ref1}
	Let $\mathcal{N}:=\mathop{\oplus}\limits_{n\geqslant 0}\mathbb{K}Q^{n} $, with the basis element corresponding to
	$(P,C)\in Q^{n}$ denoted by $H_{(P,C)}$ when $n>0$ and 1 when $(P,C)=(\emptyset_{0},\emptyset_{0})$.\par
	The coproduct on $\mathcal{N} $ is given by$$\bigtriangleup\left(H_{(P,C)}\right) =1\otimes H_{(P,C)}+H_{(P,C)}\otimes 1+\mathop{\sum}\limits_{K:rankK=1, PK\leq C}H_{(P_{1},C_{1})}\otimes H_{(P_{2},C_{2})}, $$where $b_{K}(KP)=P_{1}\times P_{2}$, $b_{K}(KC)=C_{1}\times C_{2}. $\par
	The product on $\mathcal{N} $  is given by
	$$H_{(P_{1},C_{1})}\ast H_{(P_{2},C_{2})}=H_{(j_{G}(P_{1}\times P_{2})),j_{G}(C_{1}\times C_{2}))}.$$The vertex $G$ is  given by
	Axiom $(A1)$ for the composition $(g_{1}, g_{2})$, where $P_{i} \in \Sigma^{g_{i}}$, $i=1,2$.
	\item[$\bullet$]\cite[Sections 6.5.7 and 6.8.6]{ref1}
	Let $\mathcal{S}:=\mathop{\oplus}\limits_{n\geqslant 0}\mathbb{K}\left( \mathcal{C}^{n}\times\mathcal{C}^{n}\right) ^{\ast}$, with the basis element corresponding to
	$(C,D)\in \mathcal{C}^{n}\times\mathcal{C}^{n}$ denoted by $F_{(C,D)}$ when $n>0$ and 1 when $(C,D)=(\emptyset_{0},\emptyset_{0})$.\par
	The coproduct on $\mathcal{S} $ is given by$$\bigtriangleup\left(F_{(C,D)}\right) =1\otimes F_{(C,D)}+F_{(C,D)}\otimes 1+\mathop{\sum}\limits_{K:rankK=1, K\leq D}F_{(C_{1},D_{1})}\otimes F_{(C_{2},D_{2})}, $$where $b_{K}(KC)=C_{1}\times C_{2}$, $ b_{K}(D)=D_{1}\times D_{2}. $\par
The product on $\mathcal{S} $ is given by
	$$F_{(C_{1},D_{1})}\ast F_{(C_{2},D_{2})}=\mathop{\sum}\limits_{D:GD=j_{G}(D_{1}\times D_{2})}F_{((j_{G}(C_{1}\times C_{2})),D)}.$$The vertex $G$ is  given by
	Axiom $(A1)$ for the composition $(g_{1}, g_{2})$, where $C_{i} \in \Sigma^{g_{i}}$, $i=1,2$.
	\item[$\bullet$]\cite[Sections 6.5.8 and 6.8.7]{ref1}
	Let $\mathcal{R}:=\mathop{\oplus}\limits_{n\geqslant 0}\mathbb{K} \left( \mathcal{C}^{n}\times\mathcal{C}^{n}\right) $, with the basis element corresponding to
	$(D,C)\in \mathcal{C}^{n}\times\mathcal{C}^{n}$ denoted by $K_{(D,C)}$ when $n>0$ and 1 when $(D,C)=(\emptyset_{0},\emptyset_{0})$.\par
	The coproduct on $\mathcal{R} $ is given by$$\bigtriangleup\left(K_{(D,C)}\right) =1\otimes K_{(D,C)}+K_{(D,C)}\otimes 1+\mathop{\sum}\limits_{K:rankK=1, K\leq D}K_{(D_{1},C_{1})}\otimes K_{(D_{2},C_{2})}, $$where $b_{K}(KC)=C_{1}\times C_{2}$, $ b_{K}(D)=D_{1}\times D_{2}. $\par
	The product on $\mathcal{R} $ is given by
	$$K_{(D_{1},C_{1})}\ast K_{(D_{2},C_{2})}=\mathop{\sum}\limits_{D:GD=j_{G}(D_{1}\times D_{2})}K_{(D,j_{G}(C_{1}\times C_{2}))}.$$The vertex $G$ is  given by
	Axiom $(A1)$ for the composition $(g_{1}, g_{2})$,  where $C_{i} \in \Sigma^{g_{i}}$, $i=1,2$.
	\item[$\bullet$]\cite[Sections 5.6.5-5.6.8, 6.5.10 and 6.8.9]{ref1}
	Let $\mathcal{A}_{\mathcal{L}}=\mathop{\oplus}\limits_{n\geqslant 0}\mathbb{K} L^{n}$, $\mathcal{A}_{\mathcal{Z}}=\mathop{\oplus}\limits_{n\geqslant 0}\mathbb{K} Z^{n}$, $\mathcal{A}_{\mathcal{L}^{\ast}}=\mathop{\oplus}\limits_{n\geqslant 0}\mathbb{K} (L^{n})^{\ast}$, $\mathcal{A}_{\mathcal{Z}^{\ast}}=\mathop{\oplus}\limits_{n\geqslant 0}\mathbb{K} (Z^{n})^{\ast}$, with the basis element denoted by $h_{X}$, $h_{L}$, $m_{X}$, $m_{L}$ respectively, for $X\in L^{n}$, $L\in Z^{n}$ when $n>0$ and 1 otherwise.\par
	Define following maps:
	\begin{equation*}
		\begin{split}
			supp^{\ast}: \mathcal{A}_{\mathcal{L}^{\ast}}&\rightarrow\mathcal{P},\\
			m_{X}&\mapsto \mathop{\sum}\limits_{F: supp F=X}M_{F};\\
			lune^{\ast}: \mathcal{A}_{\mathcal{Z}^{\ast}}&\rightarrow\mathcal{Q},\\m_{L}&\mapsto \mathop{\sum}\limits_{(F,D): lune (F,D)=L}M_{(F,D)};\\
			supp:\mathcal{M}&\rightarrow\mathcal{A}_{\mathcal{L}},\\H_{P}&\mapsto h_{supp P};\\
			lune:\mathcal{N}&\rightarrow\mathcal{A}_{\mathcal{Z}},\\H_{(P,C)}&\mapsto h_{lune(P,C)}.\\
		\end{split}
	\end{equation*}
\end{itemize}
\par For every graded vector space in the form of $\mathop{\oplus}\limits_{n\geqslant 0}\mathbb{K}\left( V^{n}\right)$ defined in this subsection, define a counit $\varepsilon $ by $\varepsilon(1)=1$, $ \varepsilon(\mathop{\oplus}\limits_{n\geqslant 1}\mathbb{K}\left( V^{n}\right))=\left\lbrace 0\right\rbrace $. In addition, we set 1 to be the unit of the product $\ast$ and to satisfy $\bigtriangleup(1)=1\otimes1$. Then the unit map is a coalgebra morphism.\par
From definitions, it's easy to see $lune^{\ast} $ and $supp^{\ast} $ are injective while $lune $ and  $supp $ are surjective.
The coproducts and products on $\mathcal{A}_{\mathcal{L}}$, $\mathcal{A}_{\mathcal{Z}}$, $\mathcal{A}_{\mathcal{L}^{\ast}}$ and $\mathcal{A}_{\mathcal{Z}^{\ast}}$ are also given in \cite[Chapter 6]{ref1}, but with the following theorem one can see $\mathcal{A}_{\mathcal{L}}$ and $\mathcal{A}_{\mathcal{Z}}$ as quotient algebras and quotient coalgebras of $\mathcal{M}$ and $\mathcal{N}$  respectively, and see $\mathcal{A}_{\mathcal{L}^{\ast}}$ and $\mathcal{A}_{\mathcal{Z}^{\ast}}$ as  subalgebras and subcoalgebras of $\mathcal{P}$ and  $\mathcal{Q}$ respectively.
\begin{Theorem}\label{main}\cite[Theorems 6.1.1 and 6.1.2]{ref1}
	Let $\left\lbrace \Sigma^{n}  \right\rbrace_{n\geq0}$ be a family of LRBs that satisfies all coalgebra
	axioms $(C1)-(CP)$ and all algebra
	axioms $(A1)-(AP)$. Then the following diagram is a diagram of coalgebras and algebras.
	\begin{equation}\label{taget}
		\begin{tikzcd}
			\mathcal{M} \arrow[r] \arrow[dd, "supp"', two heads] \arrow[rrdd] & \mathcal{N} \arrow[rr] \arrow[rrrdd] \arrow[dd, "lune", two heads] &                                                                                                         & \mathcal{R} \arrow[rd] &                            \\
			&                                                                                      &                                                                                                         &                             & \mathcal{S} \arrow[d]        \\
			\mathcal{A}_{\mathcal{L}} \arrow[r] \arrow[rrd]                              & \mathcal{A}_{\mathcal{Z}} \arrow[rrrd]                       & \mathcal{A}_{\mathcal{Z}^{\ast}} \arrow[rr, "lune^{\ast}" description, hook] \arrow[d] &                             & \mathcal{Q} \arrow[d] \\
			&                                                                                      & \mathcal{A}_{\mathcal{L}^{\ast}} \arrow[rr, "supp^{\ast}"', hook]                                       &                             & \mathcal{P}               
		\end{tikzcd}
	\end{equation}
	(Readers can refer to \cite[Chapter 6]{ref1} to know the definitions of other maps. We only need to know they are maps of algebras and coalgebras and that they can make this diagram commute.)
\end{Theorem}
In light of the above theorem, Aguiar and Mahajan asked a question \cite[Section 6.1]{ref1} : ``For a family $\left\lbrace \Sigma^{n}  \right\rbrace_{n\geq 0}$ of LRBs, can one give compatibility axioms such that if a family satisfies the coalgebra, algebra and compatibility axioms then Diagram \ref{taget} is a diagram of Hopf algebras?'' This paper is to solve this problem.
	\section{Compatibility axioms}	
		\subsection{Task}
		Now we consider our task. The following lemma will be useful when we deal with the maps.
		\begin{Lemma}\cite[Lemma 4.0.4]{ref3}\label{map}
			Let $H, L$ be Hopf algebras with antipodes $S_{H}$, $S_{L}$. If $V: H\rightarrow L$ is a bialgebra map, then it is a Hopf morphism (that is $f\circ S_{H}=S_{L}\circ f$).
		\end{Lemma}
		Our goal is to construct additional axioms to make Diagram \ref{taget} a diagram of Hopf algebras. With the following lemma, we can simplify our task a lot.
		\begin{Lemma}\label{task}
			If $\mathcal{M}$, $\mathcal{N}$,  $\mathcal{R}$, $\mathcal{S}$, $\mathcal{Q}$ and $\mathcal{P}$ are bialgebras, then Diagram \ref{taget} is a diagram of Hopf algebras.
		\end{Lemma}
		\begin{Proof}
			 Due to \cite[Lemma 14]{ref5}, for any graded connected bialgebra, one can get a formula for the antipode to make it a Hopf algebra.
			From axiom $(C1)$, axiom $(A1)$ and the definitions of the products and coproducts of $\mathcal{M}$, $\mathcal{N}$,  $\mathcal{R}$, $\mathcal{S}$, $\mathcal{Q}$, and $\mathcal{P}$, we can find they are all connected graded bialgebras. Accordingly they are all Hopf algebras. \par $\mathcal{A}_{\mathcal{L}^{\ast}}$, $\mathcal{A}_{\mathcal{Z}^{\ast}}$ inherit  algebra and  coalgebra structures from $\mathcal{P}$,  $\mathcal{Q}$ respectively, so they also satisfy the compatibility condition. \par
			$\mathcal{A}_{\mathcal{L}}=\mathcal{M}/ker(supp)$. $ker(supp)$ is both an ideal and a coideal of $\mathcal{M} $, so it is a biideal, which means $\mathcal{A}_{\mathcal{L}} $  is a bialgebra. $\mathcal{A}_{\mathcal{Z}} $ is bialgebra similarly.\par
			From definitions it's easy to see $lune^{\ast} $, $supp^{\ast} $, $lune$ and $supp$ are graded, so $\mathcal{A}_{\mathcal{L}}$, $\mathcal{A}_{\mathcal{Z}}$, $\mathcal{A}_{\mathcal{L}^{\ast}}$, $\mathcal{A}_{\mathcal{Z}^{\ast}}$ are also graded connected bialgebras, which means they are all Hopf algebras.\par From Lemma \ref{map}, Theorem \ref{main}  and the definition of bialgebra morphism we have that all the maps in Diagram \ref{taget} are Hopf algebra maps. $\qedsymbol$
		\end{Proof}
		Therefore our task is to give compatibility axioms to make $\mathcal{M}$, $\mathcal{N}$,  $\mathcal{R}$, $\mathcal{S}$, $\mathcal{Q}$ and $\mathcal{P}$ bialgebras.
		\subsection{Notation}
		Before showing compatibility axioms, we  define some notations for convenience.
For all $n\geq 0$, let: \begin{equation*}
	\begin{split}
		b^{n}:\Sigma^{n}&\rightarrow\Sigma^{0}\times\Sigma^{n},\\
		F&\mapsto\emptyset_{0}\times F,\\
		B^{n}:\Sigma^{n}&\rightarrow\Sigma^{n}\times\Sigma^{0},\\
		F&\mapsto F\times\emptyset_{0}.\\
	\end{split}
\end{equation*}	
\par Let $T^{n}=\left\lbrace b_{K}|K\in \Sigma^{n},rank K=1\right\rbrace \cup\left\lbrace b^{n}, B^{n}\right\rbrace $, $$T^{n}_{F}=\left\lbrace b_{K}|K\in \Sigma^{n},rank K=1, K\leq F\right\rbrace \cup\left\lbrace b^{n}, B^{n}\right\rbrace, $$ for $n\geq 1$, $F \in \Sigma^{n}$. Let $T^{0}=T^{0}_{\emptyset_{0}}=\left\lbrace b^{0}\right\rbrace .$
\par For  $f\in T^{n}$,  $P \in \Sigma^{n}$, let 	\begin{equation*}
	\hat{f}(P)=\left\{\begin{aligned}
		b_{K}(KP),\quad &if\quad f=b_{K}\quad for\quad some \quad K,\\
		f(P),\quad &otherwise.
	\end{aligned}
	\right.
\end{equation*}
\par Let \begin{equation*}
	K_{f}=\left\{\begin{aligned}
		K,\quad &if\quad f=b_{K}\quad for\quad some \quad K,\\
		\emptyset_{n},\quad &otherwise.
	\end{aligned}
	\right.
\end{equation*}
Then $ \hat{f}(P)=f(K_{f}P)$, and$$f\in T^{n}_{F}\Leftrightarrow F\in \Sigma^{n}_{K_{f}
}\Leftrightarrow K_{f}\leq F. $$
\subsection{Axioms}
Now we give our compatibility axioms. Considering that we want to use them to make bialgebras, we denoted them by $B_{i}$. \hspace*{\fill} \\
\textbf{Axiom }($B$1).\label{b1} For any $n, m\geq 0$, there exists a bijection $g_{n,m}:T^{n}\times T^{m}\rightarrow T^{n+m}$. If $f:\Sigma^{n}_{K_{f}
}\rightarrow\Sigma^{n_{1}}\times\Sigma^{n_{2}},$ $f':\Sigma^{n}_{K_{f'}
}\rightarrow\Sigma^{m_{1}}\times\Sigma^{m_{2}}$, let $f''=g_{n,m}(f,f')$, then $f''\in T^{n+m} $ is the unique poset isomorphism which satisfies the following two conditions:\\
(1) $Im f''=\Sigma^{n_{1}+m_{1}}\times\Sigma^{n_{2}+m_{2}},$\\
(2) $ GK_{f''}=j_{G}(K_{f}\times K_{f'})$, where $G$ is given by Axiom $(A1)$ for the composition $(n,m)$.\\ \hspace*{\fill} \\
\textbf{Axiom }($B$2).\label{b3} For all $n, m\geq 0$,  $f\in T^{n},$ $f'\in T^{m}$, 
let $f''=g_{n,m}(f,f')\in T^{n+m}$.  $f''$ satisfies that
$$\hat{f''}(j_{G}(P\times P'))=j_{G_{1}}(P_{1}\times P'_{1})\times j_{G_{2}}(P_{2}\times P'_{2}),$$
for any $P\in \Sigma^{n}$, $P'\in \Sigma^{m}$, where $\hat{f}(P)=P_{1}\times P_{2}\in\Sigma^{n_{1}}\times\Sigma^{n_{2}},$ 
and $\hat{f'}(P')=P'_{1}\times P'_{2}\in\Sigma^{m_{1}}\times\Sigma^{m_{2}}.$ The vertex 
 $G$ is given by Axiom ($A$1) for the composition $(n,m)$, and $G_{i} $ is given by Axiom ($A$1) for the composition $(n_{i},m_{i})$, $i=1,2$. In other words: $$\hat{f''}\circ j_{G}=(j_{G_{1}}\times j_{G_{2}})\circ (id_{n_{1}}\times\tau_{n_{2},m_{1}}\times id_{m_{2}})\circ(\hat{f}\times\hat{f'}),$$ where $ \tau_{n_{2},m_{1}}$ is the twist map $\Sigma^{n_{2}}\times\Sigma^{m_{1}}\rightarrow\Sigma^{m_{1}}\times\Sigma^{n_{2}} $ given by $(x,y)\mapsto (y,x)$.
\\\\ \hspace*{\fill} 
\subsection{Proof}
Now we will prove that our compatibility axioms $(B1)-(B2)$ together with the coalgebra axioms $(C1)-(CP)$ and algebra axioms $(A1)-(AP)$ can make Diagram \ref{taget} is a diagram of Hopf algebras.
Let $\left\lbrace \Sigma^{n}  \right\rbrace_{n\geq0}$ be a family of LRBs that satisfies all
axioms $(C1)-(CP)$, $(A1)-(AP)$ and $(B1)-(B2)$. With Lemma \ref{task}, we only need to prove  $\mathcal{M}$, $\mathcal{N}$,  $\mathcal{R}$, $\mathcal{S}$, $\mathcal{Q}$ and $\mathcal{P}$ are bialgebras. Before turning to the main results, we first prove two lemmas. The following lemmas can be seen from the compatibility axioms.
\begin{Lemma}
	For all $n, m\geq 0$,  $f\in T^{n},$ $f'\in T^{m}$,\\
	let $f''=g_{n,m}(f,f')\in T^{n+m}$. If $Imf=\Sigma^{n_{1}}\times\Sigma^{n_{2}}$, $Imf'=\Sigma^{m_{1}}\times\Sigma^{m_{2}}$,  then $$\hat{f''}(G)=G_{1}\times G_{2}.$$The vertex  $G$ is given by Axiom ($A$1) for the composition $(n,m)$, where $G_{i} $ is given by Axiom $(A1)$ for the composition $(n_{i},m_{i})$.\label{G}
\end{Lemma}
\begin{Proof}
	Apply axiom $(B2)$ to $\emptyset_{n}\times \emptyset_{m}$.\par $f:\Sigma^{n}_{K_{f}
	}\rightarrow\Sigma^{n_{1}}\times\Sigma^{n_{2}}$ is a poset isomorphism and $K_{f}\emptyset_{n}=K_{f}$, so
$$	\hat{f}(\emptyset_{n})=f(K_{f}\emptyset_{n})=f(K_{f})=\emptyset_{n_{1}}\times \emptyset_{n_{2}}.$$ \par Similarly, we have $\hat{f'}(\emptyset_{m})=\emptyset_{m_{1}}\times \emptyset_{m_{2}} $. \par Since $j_{G}:\Sigma^{n}\times\Sigma^{m}\rightarrow\Sigma^{n+m}_{G}$ is a poset isomorphism, we have $$j_{G}(\emptyset_{n}\times \emptyset_{m})=G.$$\par Similarly, $j_{G_{i}}(\emptyset_{n_{i}}\times \emptyset_{m_{i}})=G_{i} $ for $i=1,2$. Therefore $$ \hat{f''}(G)=\hat{f''}(j_{G}(\emptyset_{n}\times \emptyset_{m}))=j_{G_{1}}(\emptyset_{n_{1}}\times \emptyset_{m_{1}})\times j_{G_{2}}(\emptyset_{n_{2}}\times \emptyset_{m_{2}})=G_{1}\times G_{2}.$$ $\qedsymbol$
	\end{Proof}

\begin{Lemma}\label{b2}
	For all $n, m\geq 0$, $F\in \Sigma^{n}$, $F'\in \Sigma^{m}$,
there exists a bijection $\varpi$ from $$L=\left\lbrace (f,f',\widetilde{F_{1}},\widetilde{F_{2}})|f\in T^{n}_{F}, f'\in T^{m}_{F'},  G_{i}\widetilde{F_{i}}=j_{G_{i}}(F_{i}\times F'_{i}), i=1,2\right\rbrace, $$where $f(F)=F_{1}\times F_{2}\in\Sigma^{n_{1}}\times\Sigma^{n_{2}}, f'(F')=F'_{1}\times F'_{2}\in\Sigma^{m_{1}}\times\Sigma^{m_{2}},$ and $G_{i}$ is given by Axiom $(A1)$ for the composition $(n_{i},m_{i})$,\\
 to$$R=\left\lbrace(\widetilde{F},f'')|G\widetilde{F}=j_{G}(F\times F'), f''\in T^{n+m}_{\widetilde{F}}\right\rbrace,$$ where $ G$ is given by Axiom $(A1)$ for the composition $(n,m)$.
\end{Lemma} 
\begin{Proof}
	$\varpi$ is given by$$ (f,f',\widetilde{F_{1}},\widetilde{F_{2}})\mapsto (\widetilde{F},f''), $$where 
	 $f''=g_{n,m}(f,f')$, $\widetilde{F}=(f'')^{-1}(\widetilde{F_{1}}\times\widetilde{F_{2}})$. First we will prove $(\widetilde{F},f'')\in R. $ 
	 \begin{equation*}
			\begin{split}
				f''(K_{f''}G\widetilde{F})&=f''(K_{f''}G)f''(\widetilde{F})\\&=(G_{1}\times G_{2})(\widetilde{F_{1}}\times\widetilde{F_{2}})\\&=(G_{1}\widetilde{F_{1}})\times (G_{2}\widetilde{F_{2}})\\&=j_{G_{1}}(F_{1}\times F'_{1})\times j_{G_{2}}(F_{2}\times F'_{2})\\&=\hat{f''}(j_{G}(F\times F'))\\&=f''(K_{f''}j_{G}(F\times F')).
			\end{split}
		\end{equation*}
		The first equality follows from axiom $(CP)$ and the definition of $T^{n+m}$. The second one follows from Lemma \ref{G}. The fifth one uses axiom $(B2)$.\par
		Since $f''$ is injective, we have $$K_{f''}G\widetilde{F}=K_{f''}j_{G}(F\times F'). $$\par Premultiplying by G, then we have \begin{equation}\label{key}
		 GK_{f''}G\widetilde{F}=GK_{f''}j_{G}(F\times F').
		\end{equation} The left side of \ref{key} is $$GK_{f''}G\widetilde{F}=GK_{f''}\widetilde{F}=G\widetilde{F}$$since $\Sigma^{n+m}$ is a LRB and $K_{f''}\leq\widetilde{F}.$ The right side of \ref{key} is 
		\begin{equation*}
			\begin{split}
				GK_{f''}j_{G}(F\times F')&=j_{G}(K_{f}\times K_{f'})j_{G}(F\times F')\\&=j_{G}(K_{f}F\times K_{f'}F')\\&=j_{G}(F\times F').
			\end{split}
		\end{equation*}
		The first equality follows from axiom $(B1)$, while the second one is a consequence of axiom $(AP)$, and the third one follows that $K_{f}\leq F$ and $K_{f'}\leq F'.$ Then we have $$G\widetilde{F}=j_{G}(F\times F'),$$ that is, $(\widetilde{F},f'')\in R.$ The map $\varpi :L\rightarrow R$ given above is injective since $g_{n,m}$ is bijective. 
\par For any pair $(\widetilde{F},f'')\in R,$ let $f''(\widetilde{F})=\widetilde{F_{1}}\times\widetilde{F_{2}}, g_{n,m}^{-1}=(f,f').$ Then\begin{equation*}
		\begin{split}
			(G_{1}\widetilde{F_{1}})\times (G_{2}\widetilde{F_{2}})&=(G_{1}\times G_{2})(\widetilde{F_{1}}\times\widetilde{F_{2}})\\&=f''(K_{f''}G)f''(\widetilde{F})\\&=f''(K_{f''}G\widetilde{F})\\&=f''(K_{f''}j_{G}(F\times F'))\\&=\hat{f''}(j_{G}(F\times F'))\\&=j_{G_{1}}(F_{1}\times F'_{1})\times j_{G_{2}}(F_{2}\times F'_{2}).
		\end{split}
	\end{equation*}
	The second equality follows from Lemma \ref{G}. The third one uses axiom $(CP)$ and the definition of $T^{n+m}$.
	For the sixth equality, we apply axiom $(B2)$ to $F\times F'$ .
	Therefore $(f,f',\widetilde{F_{1}},\widetilde{F_{2}}) \in L.$
	Thus we have proved $\varpi$ is surjective. $\qedsymbol$
\end{Proof}
In what follows, we will show the compatibility of the products and the coproducts of $\mathcal{P}$, $\mathcal{M}$, $\mathcal{Q}$, $\mathcal{N}$,  $\mathcal{R}$ and  $\mathcal{S}$ in turn, that is, they are bialgebras. A direct result of the above lemma is the following theorem.
\begin{Theorem}\label{p}
	$(\mathcal{P},\bigtriangleup,\ast)$ is a bialgebra.
\end{Theorem}
\begin{Proof}
	We only need to prove the compatibility of $\ast$ and $\bigtriangleup$. Let $M_{\emptyset_{0}}=1.$ Rewrite $\bigtriangleup$ : For $F\in \Sigma^{n},$ $$\bigtriangleup(M_{F})=\mathop{\sum}\limits_{f\in T^{n}_{F}}M_{F_{1}}\otimes M_{F_{2}}, $$where $f(F)=F_{1}\times F_{2}$.\par Then for all $n, m\geq 0$, $F\in \Sigma^{n}$, $F'\in \Sigma^{n}$,
	\begin{equation*}
		\begin{split}
			\bigtriangleup(M_{F})\ast\bigtriangleup(M_{F'})&=(\mathop{\sum}\limits_{f\in T^{n}_{F}}M_{F_{1}}\otimes M_{F_{2}})\ast(\mathop{\sum}\limits_{f'\in T^{m}_{F'}}M_{F'_{1}}\otimes M_{F'_{2}})\\
			&=\mathop{\sum}\limits_{f\in T^{n}_{F},f'\in T^{m}_{F'}}(M_{F_{1}}\ast M_{F'_{1}})\otimes( M_{F_{2}},M_{F,_{2}})\\
		&=\mathop{\sum}\limits_{\substack{f\in T^{n}_{F},f'\in T^{m}_{F'},\\\widetilde{F_{i}}:G_{i}\widetilde{F_{i}}=j_{G_{i}}(F_{i}\times F'_{i}),\\i=1,2}}	M_{\widetilde{F_{1}}}\otimes M_{\widetilde{F_{2}}}, \\
		\end{split}
	\end{equation*}
	where $f(F)=F_{1}\times F_{2}\in\Sigma^{n_{1}}\times\Sigma^{n_{2}},$ and  $f'(F')=F'_{1}\times F'_{2}\in\Sigma^{m_{1}}\times\Sigma^{m_{2}}$. 
 $G_{i} $ is given by Axiom $(A1)$ for the composition $(n_{i},m_{i})$.\\
\begin{equation*}
	\begin{split}\bigtriangleup(M_{F}\ast M_{F'})&=\bigtriangleup(\mathop{\sum}\limits_{\widetilde{F}:G\widetilde{F}=j_{G}(F\times F')}M_{\widetilde{F}})\\
		&=\mathop{\sum}\limits_{\substack{\widetilde{F}:G\widetilde{F}=j_{G}(F\times F'),\\f''\in T^{n+m}_{\widetilde{F}}}}M_{\widetilde{F_{1}}}\otimes M_{\widetilde{F_{2}}}, 
	\end{split}
\end{equation*}
where $f''(\widetilde{F})=\widetilde{F}_{1}\times \widetilde{F}_{2}$, and $G$ is given by Axiom $(A1)$ for the composition $(n,m)$.\par From 
axiom $(B1)$ and Lemma \ref{b2} we have $$\bigtriangleup(M_{F})\ast\bigtriangleup(M_{F'})= \bigtriangleup(M_{F}\ast M_{F'}).$$ $\qedsymbol$
\end{Proof}
As for $ \mathcal{M}$, we mainly need Axiom ($B$2).
\begin{Theorem}
	$(\mathcal{M},\bigtriangleup,\ast)$ is a bialgebra.
\end{Theorem}
\begin{Proof}
		We only need to prove the compatibility of $\ast$ and $\bigtriangleup$. Let $H_{\emptyset_{0}}=1.$ Rewrite $\bigtriangleup$ : For $P\in \Sigma^{n},$ $$\bigtriangleup(H_{P})=\mathop{\sum}\limits_{f\in T^{n}}H_{P_{1}}\otimes H_{P_{2}}, $$where $\hat{f}(P)=P_{1}\times P_{2}$.\par Then for all $n, m\geq 0$, $P\in \Sigma^{n}$, $P'\in \Sigma^{m}$,
\begin{equation*}
	\begin{split}
		\bigtriangleup(H_{P})\ast\bigtriangleup(H_{P'})&=(\mathop{\sum}\limits_{f\in T^{n}}H_{P_{1}}\otimes H_{P_{2}})\ast(\mathop{\sum}\limits_{f'\in T^{m}}H_{P'_{1}}\otimes H_{P'_{2}})\\
		&=\mathop{\sum}\limits_{f\in T^{n},f'\in T^{m}}(H_{P_{1}}\ast H_{P'_{1}})\otimes( H_{P_{2}},H_{P,_{2}})\\
		&=\mathop{\sum}\limits_{\substack{f\in T^{n},f'\in T^{m},\\P_{i}:G_{i}\widetilde{P_{i}}=,i=1,2}}	H_{j_{G_{1}}(P_{1}\times P'_{1})}\otimes H_{j_{G_{2}}(P_{2}\times P'_{2})}, \\
	\end{split}
\end{equation*}
where $\hat{f}(P)=P_{1}\times P_{2}\in\Sigma^{n_{1}}\times\Sigma^{n_{2}}$, and  $\hat{f'}(P')=P'_{1}\times P'_{2}\in\Sigma^{m_{1}}\times\Sigma^{m_{2}}.$ 
$G_{i} $ is given by Axiom ($A$1) for the composition $(n_{i},m_{i})$.\\
\begin{equation*}
	\begin{split}\bigtriangleup(H_{P}\ast H_{P'})&=\bigtriangleup H_{j_{G}(P\times P')}\\
		&=\mathop{\sum}\limits_{f''\in T^{n+m}}H_{\widetilde{P_{1}}}\otimes H_{\widetilde{P_{2}}}, 
	\end{split}
\end{equation*}
where $\hat{f''}(j_{G}(P\times P'))=\widetilde{P}_{1}\times \widetilde{P}_{2}$, and $G$ is given by Axiom ($A$1) for the composition $(n,m)$.\par From 
Axiom ($B$1) and Axiom ($B$2) we have $$\bigtriangleup(H_{P})\ast\bigtriangleup(H_{P'})= \bigtriangleup(H_{P}\ast H_{P'}).$$  $\qedsymbol$
\end{Proof}
%
The coalgebra structure and the algebra structure of $\mathcal{Q} $ in the first coordinate are similar to those of $\mathcal{P} $, but we need to show the correspondence of the second coordinate.
\begin{Theorem}
	$(\mathcal{Q},\bigtriangleup,\ast)$ is a bialgebra.
\end{Theorem}
\begin{Proof}
	We only need to prove the compatibility of $\ast$ and $\bigtriangleup$. Let $M_{(\emptyset_{0},\emptyset_{0})}=1.$ Rewrite $\bigtriangleup$ : For $(F,D)\in Q^{n},$ $$\bigtriangleup(M_{(F,D)})=\mathop{\sum}\limits_{f\in T^{n}_{F}}M_{(F_{1},D_{1})}\otimes M_{(F_{2},D_{2})}, $$where $f(F)=F_{1}\times F_{2}$, and $f(D)=D_{1}\times D_{2}$.
	\par Then for all $n, m\geq 0$, $(F,D)\in Q^{n}$, $(F',D')\in Q^{m}$,
	\begin{equation*}
		\begin{split}
			&\bigtriangleup(M_{(F,D)})\ast\bigtriangleup(M_{(F',D')})\\=&(\mathop{\sum}\limits_{f\in T^{n}_{F}}M_{(F_{1},D_{1})}\otimes M_{(F_{2},D_{2})})\ast(\mathop{\sum}\limits_{f'\in T^{m}_{F'}}M_{(F'_{1},D'_{1})}\otimes M_{(F'_{2},D'_{2})})\\
			=&\mathop{\sum}\limits_{f\in T^{n}_{F},f'\in T^{m}_{F'}}(M_{(F_{1},D_{1})}\ast M_{(F'_{1},D'_{1})})\otimes(M_{(F_{2},D_{2})}\ast M_{(F'_{2},D'_{2})})\\
			=&\mathop{\sum}\limits_{\substack{f\in T^{n}_{F},f'\in T^{m}_{F'},\\\widetilde{F_{i}}:G_{i}\widetilde{F_{i}}=j_{G_{i}}(F_{i}\times F'_{i}),\\i=1,2}}	M_{(\widetilde{F_{1}},\widetilde{F_{1}}j_{G_{1}}(D_{1}\times D'_{1}))}\otimes M_{(\widetilde{F_{2}},\widetilde{F_{2}}j_{G_{2}}(D_{2}\times D'_{2}))}, \\
		\end{split}
	\end{equation*}
	where $f(F)=F_{1}\times F_{2}\in\Sigma^{n_{1}}\times\Sigma^{n_{2}},$ $f'(F')=F'_{1}\times F'_{2}\in\Sigma^{m_{1}}\times\Sigma^{m_{2}},$ $f(D)=D_{1}\times D_{2},$ and $f'(D')=D'_{1}\times D'_{2}.$
	$G_{i} $ is given by Axiom ($A$1) for the composition $(n_{i},m_{i})$.\\
		\begin{equation*}
		\begin{split}\bigtriangleup(M_{(F,D)}\ast M_{(F',D')})&=\bigtriangleup(\mathop{\sum}\limits_{\widetilde{F}:G\widetilde{F}=j_{G}(F\times F')}M_{(\widetilde{F},\widetilde{F}j_{G}(D\times D'))})\\
			&=\mathop{\sum}\limits_{\substack{\widetilde{F}:G\widetilde{F}=j_{G}(F\times F'),\\f''\in T^{n+m}_{\widetilde{F}}}}M_{(\widetilde{F_{1}},\widetilde{D_{1}})}\otimes M_{(\widetilde{F_{2}},\widetilde{D_{2}})}, \\
		\end{split}
	\end{equation*}
	where $f''(\widetilde{F})=\widetilde{F}_{1}\times \widetilde{F}_{2}$, and  $f''(\widetilde{F}j_{G}(D\times D'))=\widetilde{D}_{1}\times \widetilde{D}_{2}$. $G$ is given by Axiom ($A$1) for the composition $(n,m)$.\par
	We claim  $$f''(\widetilde{F}j_{G}(D\times D'))=(\widetilde{F_{1}}j_{G_{1}}(D_{1}\times D'_{1}))\times(\widetilde{F_{2}}j_{G_{2}}(D_{2}\times D'_{2})),$$ where $f''=g_{n,m}(f,f')$. With this claim and Lemma \ref{b2},\\
	we have $$ \bigtriangleup(M_{(F,D)})\ast\bigtriangleup(M_{(F',D')})=\bigtriangleup(M_{(F,D)}\ast M_{(F',D')}).$$\par Now we prove the claim:
	\begin{equation*}
		\begin{split}
			f''(\widetilde{F}j_{G}(D\times D'))&=f''(K_{f''}\widetilde{F}j_{G}(D\times D'))\\&=f''(K_{f''}\widetilde{F}K_{f''}j_{G}(D\times D'))\\&=f''(K_{f''}\widetilde{F})f''(K_{f''}j_{G}(D\times D'))\\&=f''(\widetilde{F})\hat{f''}(j_{G}(D\times D'))\\&=(\widetilde{F}_{1}\times \widetilde{F}_{2})(j_{G_{1}}(D_{1}\times D'_{1})\times j_{G_{2}}(D_{2}\times D'_{2}))\\&=(\widetilde{F_{1}}j_{G_{1}}(D_{1}\times D'_{1}))\times(\widetilde{F_{2}}j_{G_{2}}(D_{2}\times D'_{2})).
		\end{split}
	\end{equation*}
	The first and the fourth equalities follow from that $f''\in T^{n+m}_{\widetilde{F}} $, that is, $K_{f''}\leq\widetilde{F}. $ The second one follows from the definition of LRB. The third one follows from axiom $(CP)$ and the definition of $T^{n+m}$. For the fifth equality, we apply axiom $(B2)$ to $D\times D'$ .
 $\qedsymbol$.
\end{Proof}
The coalgebra structure and the algebra structure of $\mathcal{N} $  are similar to those of $\mathcal{M} $, but with a restriction. We need to show the correspondence of the restrictions.
\begin{Theorem}
	$(\mathcal{N},\bigtriangleup,\ast)$ is a bialgebra.
\end{Theorem}
\begin{Proof}
	We only need to prove the compatibility of $\ast$ and $\bigtriangleup$. Let $H_{(\emptyset_{0},\emptyset_{0})}=1.$ Rewrite $\bigtriangleup$ : For $(P,C)\in Q^{n}$, $$\bigtriangleup(H_{(P,C)})=\mathop{\sum}\limits_{f\in T^{n}, PK_{f}\leq C}H_{(P_{1},C_{1})}\otimes H_{(P_{2},C_{2})}, $$where $\hat{f}(P)=P_{1}\times P_{2}$ and  $\hat{f}(C)=C_{1}\times C_{2}$.\par Then for all $n, m\geq 0$, $(P,C)\in Q^{n},$  and $(P',C')\in Q^{m}$,
	\begin{equation*}
		\begin{split}
			&\bigtriangleup(H_{(P,C)})\ast\bigtriangleup(H_{(P',C')})\\=&(\mathop{\sum}\limits_{f\in T^{n},PK_{f}\leq C}H_{(P_{1},C_{1})}\otimes H_{(P_{2},C_{2})})\ast(\mathop{\sum}\limits_{f'\in T^{m},P'K_{f'}\leq C'}H_{(P'_{1},C'_{1})}\otimes H_{(P'_{2},C'_{2})})\\
			=&\mathop{\sum}\limits_{\substack{f\in T^{n},f'\in T^{m},\\PK_{f}\leq C,P'K_{f'}\leq C'}}(H_{(P_{1},C_{1})}\ast H_{(P'_{1},C'_{1})})\otimes( H_{(P_{2},C_{2})},H_{(P'_{2},C_{2})})\\
			=&\mathop{\sum}\limits_{\substack{f\in T^{n}_{P},f'\in T^{m}_{P'},\\PK_{f}\leq C,P'K_{f'}\leq C'}}	H_{(j_{G_{1}}(P_{1}\times P'_{1}),j_{G_{1}}(C_{1}\times C'_{1}))}\otimes H_{(j_{G_{2}}(P_{2}\times P'_{2}),j_{G_{2}}(C_{2}\times C'_{2}))}, \\
		\end{split}
	\end{equation*}
	where $\hat{f}(P)=P_{1}\times P_{2}\in\Sigma^{n_{1}}\times\Sigma^{n_{2}},$ $\hat{f}(C)=C_{1}\times C_{2}$, $\hat{f'}(P')=P'_{1}\times P'_{2}\in\Sigma^{m_{1}}\times\Sigma^{m_{2}},$ and $\hat{f'}(P')=P'_{1}\times P'_{2}.$
	$G_{i} $ is given by Axiom $(A1)$ for the composition $(n_{i},m_{i})$.\\
	\begin{equation*}
		\begin{split}\bigtriangleup(H_{(P,C)}\ast H_{(P',C')})&=\bigtriangleup H_{(j_{G}(P\times P'),j_{G}(C\times C'))}\\
			&=\mathop{\sum}\limits_{\substack{f''\in T^{n+m}_{j_{G}(P\times P')},\\j_{G}(P\times P')K_{f''}\leq j_{G}(C\times C')}}H_{(\widetilde{P_{1}},\widetilde{C_{1}})}\otimes H_{(\widetilde{P_{2}},\widetilde{C_{2}})}, 
		\end{split}
	\end{equation*}
	where $\hat{f''}(j_{G}(P\times P'))=\widetilde{P}_{1}\times \widetilde{P}_{2}$, and $\hat{f''}(j_{G}(C\times C'))=\widetilde{C}_{1}\times \widetilde{C}_{2}$. $G$ is given by Axiom $(A1)$ for the composition $(n,m)$.\\	We claim  $$PK_{f}\leq C,P'K_{f'}\leq C'\Leftrightarrow j_{G}(P\times P')K_{f''}\leq j_{G}(C\times C'),$$ where $f''=g_{n,m}(f,f').$\par With this claim and axiom $(B2)$,
	we have $$ \bigtriangleup(H_{(P,C)})\ast\bigtriangleup(H_{(P',C')})=\bigtriangleup(H_{(P,C)}\ast H_{(P',C')}).$$\par Now we prove the claim:\begin{equation*}
		\begin{split}
			j_{G}(P\times P')K_{f''}&=Gj_{G}(P\times P')K_{f''}\\
			&=Gj_{G}(P\times P')GK_{f''}\\&=j_{G}(P\times P')GK_{f''}\\&=j_{G}(P\times P')j_{G}(K_{f}\times K_{f'})\\&=j_{G}((PK_{f})\times (P'K_{f'}))	\end{split}		
	\end{equation*}
	The first and the third equality follows that $G\leq j_{G}(P\times P')$ . The second one follows from the definition of a LRB. The fourth equality follows from axiom $(B1)$, while the fifth one is a consequence of axiom $(AP)$. 
	  $\qedsymbol$
\end{Proof}
The coalgebra and algebra structures of $\mathcal{S} $ and $\mathcal{R} $ are combinations of those of $\mathcal{P} $ and $\mathcal{M} $, so the proof is easy.
\begin{Theorem}\label{st}
		$(\mathcal{S},\bigtriangleup,\ast)$ and $(\mathcal{R},\bigtriangleup,\ast)$ are bialgebras.
\end{Theorem}
\begin{Proof}

	We only need to prove for $(\mathcal{S},\bigtriangleup,\ast)$. Rewrite $\bigtriangleup$: For $(C,D)\in \mathcal{C}^{n}\times\mathcal{C}^{n},$
$$\bigtriangleup(F_{(C,D)})=\mathop{\sum}\limits_{f\in T^{n}_{D}}F_{(C_{1},D_{1})}\otimes F_{(C_{2},D_{2})}, $$where $\hat{f}(C)=C_{1}\times C_{2}$ and $f(D)=D_{1}\times D_{2}$. 
	\begin{equation*}
	\begin{split}
		& \bigtriangleup(F_{(C,D)})\ast\bigtriangleup(F_{(C',D')})\\=&(\mathop{\sum}\limits_{f\in T^{n}_{D}}F_{(C_{1},D_{1})}\otimes F_{(C_{2},D_{2})})\ast(\mathop{\sum}\limits_{f'\in T^{m}_{D'}}F_{(C'_{1},D'_{1})}\otimes F_{(C'_{2},D'_{2})})\\
		=&\mathop{\sum}\limits_{f\in T^{n}_{D},f'\in T^{m}_{D'}}(F_{(C_{1},D_{1})}\ast F_{(C'_{1},D'_{1})})\otimes(F_{(C_{2},D_{2})}\ast F_{(C'_{2},D'_{2})})\\
		=&\mathop{\sum}\limits_{\substack{f\in T^{n}_{D},f'\in T^{m}_{D'},\\\widetilde{D_{i}}:G_{i}\widetilde{D_{i}}=j_{G_{i}}(D_{i}\times D'_{i}),\\i=1,2}}	F_{(j_{G_{1}}(C_{1}\times C'_{1}),\widetilde{D_{1}})}\otimes F_{(j_{G_{2}}(C_{2}\times C'_{2}),\widetilde{D_{2}})}, \\
	\end{split}
\end{equation*}
where $\hat{f}(C)=C_{1}\times C_{2}\in\Sigma^{n_{1}}\times\Sigma^{n_{2}},$ $\hat{f'}(C')=C'_{1}\times C'_{2}\in\Sigma^{m_{1}}\times\Sigma^{m_{2}},$ $f(D)=D_{1}\times D_{2},$ and $f'(D')=D'_{1}\times D'_{2}.$
$G_{i} $ is given by Axiom $(A1)$ for the composition $(n_{i},m_{i})$.\\
\begin{equation*}
	\begin{split}\bigtriangleup(F_{(C,D)}\ast F_{(C',D')})&=\bigtriangleup(\mathop{\sum}\limits_{\widetilde{D}:G\widetilde{D}=j_{G}(D\times D')}F_{(j_{G}(C\times C'),\widetilde{D})})\\
		&=\mathop{\sum}\limits_{\substack{\widetilde{D}:G\widetilde{D}=j_{G}(D\times D'),\\f''\in T^{n+m}_{\widetilde{D}}}}F_{(\widetilde{C_{1}},\widetilde{D_{1}})}\otimes F_{(\widetilde{C_{2}},\widetilde{D_{2}})}, \\
	\end{split}
\end{equation*}
where $\hat{f''}(j_{G}(C\times C'))=\widetilde{C}_{1}\times \widetilde{C}_{2}$,  $f''(\widetilde{D})=\widetilde{D}_{1}\times \widetilde{D}_{2}$. $G$ is given by Axiom $(A1)$ for the composition $(n,m)$.\par
Apply Axiom ($B2)$ to $C\times C'$, we have that if $f''=g_{n,m}(f,f')$, then$$ \widetilde{C}_{1}\times \widetilde{C}_{2}=\hat{f''}(j_{G}(C\times C'))=j_{G_{1}}(C_{1}\times C'_{1})\times j_{G_{2}}(C_{2}\times C'_{2}).$$Use Lemma  \ref{b2}
we have that$$ \bigtriangleup(F_{(C,D)})\ast\bigtriangleup(M_{(F',D')})=\bigtriangleup(F_{(C,D)}\ast F_{(C',D')}).$$ $\qedsymbol$
\end{Proof}
With these results at hand, we can now embark on the main theorem.
\begin{Theorem}\label{main result}
	Let $\left\lbrace \Sigma^{n}  \right\rbrace_{n\geq0}$ be a family of LRBs that satisfies all 
	axioms $(C1)-(CP)$,($A1)-(AP$) and ($B1)-(B2$). Then Diagram \ref{taget} is a diagram of Hopf algebras.
\end{Theorem}
\begin{Proof}
	The theorem \ref{main result} is a direct consequence of Lemma \ref{task} and theorems \ref{p}-\ref{st}.
\end{Proof}
\subsection{Example}
In this section we will give an example of a family of LRBs which satisfies all the axioms $(C1)-(CP)$, ($A1)-(AP)$, and ($B1)-(B2$). The purpose of this part is to show our axioms are compatible with each other.\par
Let $n$ be a positive integer. A \emph{composition} of $[n]$ is an ordered partition $F^{1}\cdots F^{k}$ of $[n]$. In other words, $F^{1}\cdots F^{k} $ are disjoint nonempty sets
whose union is $[n]$, and their order counts. Let $ {\mathcal{B}^{n}}$ be the compositions of $[n]$. Define a product on $ {\mathcal{B}^{n}}$ to make it a LRB as follows:
If $F=F^{1}\cdots F^{k}$ and $H=H^{1}\cdots H^{m}$, then 
$$FH=(F^{1}\cap H^{1}\cdots F^{1}\cap H^{m}\cdots F^{k}\cap H^{1} \cdots F^{k}\cap H^{m})\hat{}$$
where the hat  means “deleting empty intersections”. \par

It is easy to check that ${\mathcal{B}^{n}}$ is a LRB. $F\leq H\Leftrightarrow FH=H\Leftrightarrow$ $H$ is a refinement of $F$.
Rank the compositions of $[n]$ by the numbers of ``|'', then ${\mathcal{B}^{n}}$ is a graded poset of rank $n-1$.\par
	There is a unique order-preserving map $st$ from any $n$-set
$N$ of the integers to the standard $n$-set $[n]$. Using this map, one can standardize a composition of $N$ to a composition of $[n]$. For example,
$st(9|16|57) = 5|13|24$. ``$st$'' is called the standardization map. Let $A$ be a set and $b=\lvert A\rvert$, the standard map from $A $ to $[b]$ is denoted $st_{A,b}$. \par
Let F be a set of integers, and $m$ an integer. Let $F+m:=\left\lbrace a+m|a\in F\right\rbrace $.

Now we give the definition of $b_{F}$ and $j_{F}$ for a set composition $F$.
\begin{Definition}\cite[Definition 6.3.3]{ref1}
	For any composition $F=F^{1}\cdots F^{l}$ of $[n]$, let $\mathcal{B}^{n}_{F}=\left\lbrace H\in \mathcal{B}^{n}|F\leq H\right\rbrace $. We have
	$l = deg F$. Let $f_{i}$ be the cardinality of $F_{i}$ for $1 \leq i \leq l$. The poset isomorphism
	$$b_{F}:\mathcal{B}^{n}_{F}\rightarrow \mathcal{B}^{f_{1}}\times \mathcal{B}^{f_{2}}\times \cdots\times\mathcal{B}^{f_{l}}$$
	is defined as follows.
	Let $H \in \mathcal{B}^{n}_{F}$ be a refinement of $F$. Then the image of $H$ on the $i$-th factor is
	obtained by lumping together the blocks of $H$ that refine $F_{i}$ and standardizing so that the result is a composition of the set [$f_{i}$].\par For example, for $n=8$ and  $F =14|2357|68$, the map $b_{F}$ :
	$$b_{F}:\mathcal{B}^{8}_{F}\rightarrow \mathcal{B}^{2}\times \mathcal{B}^{4}\times\mathcal{B}^{2}$$ 
	sends  $H =14|5|3|27|8|6$ to the triplet $12\times3|2|14\times 2|1$.\par
	Let ${\mathcal{B}^{0}}=\left\lbrace \emptyset_{0} \right\rbrace $, where $\left\lbrace \emptyset_{0} \right\rbrace $ is the empty composition of $\emptyset$.
\end{Definition}

\begin{Definition}\cite[Definition 6.6.2]{ref1}
	For a composition $(f_{1},\cdots,f_{k})$ of $n$, let  $F =F^{1}| \cdots |F^{k}$ be
	the composition of the set $[n]$, where $F^{1}= 12 \cdots f_{1}$, $F^{2} =f_{1} + 1 \cdots f_{1} + f_{2}$, and so
	on. Then define a poset isomorphism $$j_{F} : \mathcal{B}^{f_{1}}\times \mathcal{B}^{f_{2}}\times \cdots\times\mathcal{B}^{f_{k}} \rightarrow \mathcal{B}^{n}
	_{F} $$
	as follows.\par
	Let $F_{i}\in \mathcal{B}^{f_{i}}$ be a composition of the set [$f_{i}$]. Then the image of $F_{1}\times F_{2}\times\cdots\times F_{k}$
	is obtained by placing $F_{1} $, $F_{2}+f_{1}$, $F_{3}+(f_{1} + f_{2})$, and so on next to one another.\par
	For example, for the composition (2, 3, 3) of 8, we have $F =
	12|345|678,$ and the map $$j_{F} : \mathcal{B}^{2}\times \mathcal{B}^{3}\times \mathcal{B}^{3} \rightarrow \mathcal{B}^{8}
	_{F} $$ sends $12\times 3|2|1\times3|12$ to the
	set composition 12|5|4|3|8|67.
\end{Definition}

Aguiar and Mahajan \cite[Lemma 6.3.3 and Lemma 6.6.3]{ref1} have
told us $\left\lbrace \mathcal{B}^{n}  \right\rbrace_{n\geq 0}$ satisfies all coalgebra axioms	$(C1)-(CP)$ and algebra axioms $(A1)-(AP)$. Our main task in this part is to prove that $\left\lbrace \mathcal{B}^{n}  \right\rbrace_{n\geq 0}$ satisfies our compatibility axioms $(B1)$ and $(B2)$ as well.
 
\begin{Theorem}
	The family of $\left\lbrace \mathcal{B}^{n}  \right\rbrace_{n\geq 0}$ satisfies compatibility axioms $(B1)$ and $(B2)$.
\end{Theorem}
\begin{Proof}
	For $n\geq$0, define a map $\iota:T^{n}\rightarrow \left\lbrace S|S\subseteq [n]\right\rbrace :$
	\begin{equation*}
		\iota(f)=\left\{\begin{aligned}
			&F^{1},\quad if\quad f=b_{K}\quad for\quad some \quad K=F^{1}|F^{2},\\ &\emptyset,\quad if\quad f=b^{n},\\
			&[n],\quad if\quad f=B^{n}.
		\end{aligned}
		\right.
	\end{equation*}
	It is easy to see $\iota$ is a bijection. For any $ S\subseteq [n]$,
	\begin{equation*}
		\iota^{-1}(S)=\left\{\begin{aligned}
			&b_{S|[n]\setminus S},\quad if\quad \emptyset\subsetneq S\subsetneq [n],\\ &b^{n},\quad if\quad S=\emptyset,\\
			&B^{n},\quad if\quad S=[n].
		\end{aligned}
		\right.
	\end{equation*}Then we have $K_{f}=(\iota(f)|[n]\setminus\iota(f))\hat{}$.\par
	Define $g_{n,m}:T^{n}\times T^{m}\rightarrow T^{n+m}$ by $$(\iota^{-1}(F),\iota^{-1}(F'))\mapsto \iota^{-1}(F\cup(F'+n)),$$ where $F\subseteq[n]$, $F'\subseteq[m]$.\par  Axiom $(A1)$ for the composition $(n,m)$ gives the vertex $G=1\cdots n|n+1\cdots n+m$. Let $f=\iota^{-1}(F),$ $f'=\iota^{-1}(F')$, $f''=g_{n,m}(f,f')$, then  $K_{f}=(F|[n]\setminus F)\hat{} $, $K_{f'}=(F'|[m]\setminus F')\hat{} $, $$ K_{f''}=(F\cup(F'+n)|[n+m]\setminus(F\cup(F'+n)))\hat{},$$ $$GK_{f''}=(F|[n]\setminus F|F'+n|([m]\setminus F')+n)\hat{}=j_{G}(K_{f}\times K_{f'}).$$\par Therefore the family $\left\lbrace \mathcal{B}^{n}  \right\rbrace_{n\geq 0}$ satisfies compatibility axioms $(B1)$. 
	For any $P=P^{1}|\cdots|P^{s}\in\mathcal{B}^{n}$, $P'=P'^{1}|\cdots|P'^{t}\in\mathcal{B}^{m}$, we have$$j_{G}(P\times P')=P^{1}|\cdots|P^{s}|P'^{1}+n|\cdots|P'^{t}+n,$$
	$$K_{f}P=(F\cap P^{1}|\cdots|F\cap P^{s}|P^{1}\setminus F |\cdots|P^{s}\setminus F )\hat{}, $$	$$K_{f'}P'=(F'\cap P'^{1}|\cdots|F'\cap P'^{t}|P'^{1}\setminus F' |\cdots|(P'^{t}\setminus F' )\hat{}, $$
	\begin{equation*}
		\begin{split}
			K_{f''}j_{G}(P\times P')=&(F\cap P^{1}|\cdots|F\cap P^{s}|(F'\cap P'^{1})+n|\cdots|(F'\cap P'^{t})+n|\\&P^{1}\setminus F |\cdots|P^{s}\setminus F|(P'^{1}\setminus F')+n |\cdots|(P'^{t}\setminus F')+n )\hat{}, 
		\end{split}
	\end{equation*}
		then
	$$\hat{f}(P)=st((F\cap P^{1}|\cdots|F\cap P^{s})\hat{})\times st((P^{1}\setminus F |\cdots|P^{s}\setminus F)\hat{}),$$ $$\hat{f'}(P')=st((F'\cap P'^{1}|\cdots|F'\cap P'^{t})\hat{})\times st((P'^{1}\setminus F' |\cdots|(P'^{t}\setminus F')\hat{}),$$
	\begin{equation*}
		\begin{split}
			\hat{f''}j_{G}(P\times P')=&st((F\cap P^{1}|\cdots|F\cap P^{s}|(F'\cap P'^{1})+n|\cdots|(F'\cap P'^{t})+n)\hat{})\\&\times st((P^{1}\setminus F |\cdots|P^{s}\setminus F|(P'^{1}\setminus F')+n |\cdots|(P'^{t}\setminus F')+n )\hat{}).
		\end{split}
	\end{equation*}

	Let $n_{1}=\lvert F\rvert$, $n_{2}=n-n_{1}$, $m_{1}=\lvert F'\rvert$, $m_{2}=m-m_{1}$,  then $(n_{i},m_{i})$ gives $G_{i}=1\cdots n_{i}|n_{i}+1\cdots n_{i}+m_{i}$.
	Then consider $$
	L=j_{G_{1}}(st((F\cap P^{1}|\cdots|F\cap P^{s})\hat{})\times st((F'\cap P'^{1}|\cdots|F'\cap P'^{t})\hat{}))$$ and 
	$$R=j_{G_{2}}(st((P^{1}\setminus F |\cdots|P^{s}\setminus F)\hat{})\times st((P'^{1}\setminus F' |\cdots|(P'^{t}\setminus F')\hat{})).$$
	 $L$ gives a map  $st':F\cup(F'+n)\rightarrow[n_{1}+m_{1}]$ by
	\begin{equation*}
	st'(a)=\left\{\begin{aligned}
		&st_{F,n_{1}}(a),\quad if\quad a\in F,\\
		&st_{F',{m_{1}}}(a-n)+n_{1},\quad if\quad a\in F'+n.
	\end{aligned}
		\right.
	\end{equation*}
 $st'$ is order-preserving so $st'=st_{F\cup(F'+n),n_{1}+m_{1}}$, that is, $$L=st((F\cap P^{1}|\cdots|F\cap P^{s}|(F'\cap P'^{1})+n|\cdots|(F'\cap P'^{t})+n)\hat{}).$$\par Similarly,$$ R=st((P^{1}\setminus F |\cdots|P^{s}\setminus F|(P'^{1}\setminus F')+n |\cdots|(P'^{t}\setminus F')+n )\hat{}).$$ Then we have$$\hat{f''}j_{G}(P\times P')=L\times R. $$\par We conclude that the family $\left\lbrace \mathcal{B}^{n}  \right\rbrace_{n\geq 0}$ satisfies compatibility axioms $(B1)$. $\qedsymbol$
\end{Proof}
Apply Theorem \ref{main} to $\left\lbrace \mathcal{B}^{n}  \right\rbrace_{n\geq 0}$, then Diagram \ref{taget} specializes to  	\\
\begin{equation*}
\begin{tikzcd}
	M\Pi \arrow[r, "base^{\ast}"] \arrow[dd, "supp"', two heads] \arrow[rrdd, "\Upsilon^{\ast}"'] & N\Pi \arrow[rr, "\Theta"] \arrow[rrrdd, "\Psi"] \arrow[dd, "lune", two heads] &                                                                                                         & R\Pi \arrow[rd, "s"] &                            \\
	&                                                                                      &                                                                                                         &                             & S\Pi \arrow[d, "Road"]        \\
	\Pi_{L} \arrow[r, "base^{\ast}"] \arrow[rrd, "\Phi"']                              & \Pi_{Z} \arrow[rrrd, "\Upsilon" description]                       & \Pi_{Z^{\ast}} \arrow[rr, "lune^{\ast}" description, hook] \arrow[d, "base" description] &                             & Q\Pi \arrow[d, "base"] \\
	&                                                                                      & \Pi_{L^{\ast}} \arrow[rr, "supp^{\ast}"', hook]                                       &                             & P\Pi               
\end{tikzcd}
\end{equation*}
where the Hopf algebra $M\Pi$, $N\Pi$, etc. and the maps can be defined in a more combinatorial way \cite[Section 6.2]{ref1}. 
$P\Pi$, $ \Pi_{L^{\ast}}$ and $ \Pi_{L}$ are also studied in \cite{ref8, ref9, Montgomery1, ref33, ref45, ref68, ref72, ref74, ref83, ref100}.
The internal structure of $M\Pi$ is studied in \cite{ref75}.

\end{document}